\documentclass[leqno,12pt]{amsart}
\usepackage{a4,latexsym,amssymb,amsfonts}
\usepackage{enumerate}
\usepackage{amsthm}
\numberwithin{equation}{section}

\newtheorem{theorem}{Theorem}[section]
\newtheorem{lemma}[theorem]{Lemma}
\newtheorem{proposition}[theorem]{Proposition}
\newtheorem{corollary}[theorem]{Corollary}

\theoremstyle{definition}
\newtheorem{remark}[theorem]{Remark}

\newcommand{\rrepsp}{\frac{\rho'}{\sqrt{\nu}} }

\newcommand{\rreps}{\frac{\rho}{\sqrt{\nu}} }
\newcommand{\repn}{\frac{1}{{\nu}} }

\newcommand{\pkne}{\pi_{k,N(\nu)-k} }
\newcommand{\pllp}{\pi_{\ell,\ell'} }
\newcommand{\pnke}{\pi_{N(\nu)-k,k} }

\newcommand\dis{\displaystyle}
\newcommand\tp{\tilde{p}}

\newcommand\somma{(\ell+\ell')}

\newcommand\ql{{q_{\ell}}}
\newcommand\Ql{{Q_{\ell}}}
\newcommand\dll{d_{\ell,\ell'}}

\newcommand\hll{\cH^{\ell,\ell'}}
\newcommand\hlo{\cH^{\ell,0}}
\newcommand\zll{\ZZ_{\ell,\ell'}}
\newcommand\zlo{\ZZ_{\ell,0}}

\newcommand\pimezzi{\frac{\pi}{2}}

\newcommand\RR{{\mathbb{R}}}
\newcommand\CC{{\mathbb{C}}}
\newcommand\NN{{\mathbb{N}}}
\newcommand\ZZ{{\mathbb{Z}}}

\newcommand\TT{{\mathbb{T}}}

\newcommand\cM{{\mathcal{M}}}
\newcommand\cH{{\mathcal{H}}}
\newcommand\HH{{\mathcal{H}}}

\newcommand\cE{{\mathcal{E}}}
\newcommand\cL{{\mathcal{L}}}

\newcommand\cI{{\mathcal{I}}}

\newcommand\heise{{\it{h}}}

\newcommand\unme{\frac{1}{2}}
\newcommand\unmeko{\frac{1}{2k_{0}}}

\newcommand\unpi{\frac{1}{p}}

\newcommand\zkne{{\ZZ}_{k, N(\nu)-k}}

\def\R{{\rm I\mskip -3.5mu R}}
\def\Z{{\rm Z\mskip -13mu Z}}
\def\N{{\rm I\mskip -3.5mu N}}

\def\nmenouno{\frac{n }{2 }}
\def\tif{{\tilde{f} }}
\def\fnu{{f_{\nu} }}
\def\eeta{\underline{\bf{\eta}}}
\def\eetap{\underline{\bf{\eta'}}}
\def\thetaj{{\theta_{j} }}
\def\sfera{{S^{2n+1} }}
\newcommand\heisenunouno{{\it{h}^{n}}}
\newcommand\heiserid{{\it{h}^{n}}}
\newcommand\heisenuno{{{H}^{n}}}
\def\underz{{\bf{z}}}
\def\underw{{\bf{w}}}
\def\pmk{{P_{m,k} }}
\def\Do{{\Delta_{0} }}
\def\wnb{{\bar w}_{n+1}}
\def\wub{{\bar w}_{1}}
\def\znb{{\bar z}_{n+1}}

\def\zub{{\bar z}_{1}}
\def\zu{z_{1}}
\def\zn{z_{n+1}}

\def\uno{\bold{ 1}}

\def\zetab{\bar {{\underz}}}
\def\hllp{\cH^{\ell\ell'}}

\begin{document}

\font\fivrm=cmr5\relax


\catcode`@=11 \catcode`!=11

\expandafter\ifx\csname fiverm\endcsname\relax
  \let\fiverm\fivrm
\fi

\let\!latexendpicture=\endpicture
\let\!latexframe=\frame
\let\!latexlinethickness=\linethickness
\let\!latexmultiput=\multiput
\let\!latexput=\put

\def\@picture(#1,#2)(#3,#4){%
  \@picht #2\unitlength
  \setbox\@picbox\hbox to #1\unitlength\bgroup
  \let\endpicture=\!latexendpicture
  \let\frame=\!latexframe
  \let\linethickness=\!latexlinethickness
  \let\multiput=\!latexmultiput
  \let\put=\!latexput
  \hskip -#3\unitlength \lower #4\unitlength \hbox\bgroup}

\catcode`@=12 \catcode`!=12

\input{pictex}



\catcode`@=11 \catcode`!=11

\let\!pictexendpicture=\endpicture
\let\!pictexframe=\frame
\let\!pictexlinethickness=\linethickness
\let\!pictexmultiput=\multiput
\let\!pictexput=\put

\def\beginpicture{%
  \setbox\!picbox=\hbox\bgroup%
  \let\endpicture=\!pictexendpicture
  \let\frame=\!pictexframe
  \let\linethickness=\!pictexlinethickness
  \let\multiput=\!pictexmultiput
  \let\put=\!pictexput
  \let\input=\@@input   
primitive
  \!xleft=\maxdimen
  \!xright=-\maxdimen
  \!ybot=\maxdimen
  \!ytop=-\maxdimen}

\let\frame=\!latexframe

\let\pictexframe=\!pictexframe

\let\linethickness=\!latexlinethickness
\let\pictexlinethickness=\!pictexlinethickness

\let\\=\@normalcr
\catcode`@=12 \catcode`!=12

\vspace*{2.4in}

\title[Transferring $L^p$ eigenfunction bounds]{
Transferring $L^p$ eigenfunction bounds\\
from $\sfera$ to $\heisenunouno$}
\author{ Valentina Casarino}
\address{ Dipartimento di Matematica\\
Politecnico di Torino\\
Corso Duca degli Abruzzi 24\\10129 Torino}
\address{ 
Dipartimento di Metodi e Modelli matematici per le scienze applicate\\
Via Trieste 63, 35121 Padova}
\author{Paolo Ciatti}
\email{casarino@calvino.polito.it, ciatti@dmsa.unipd.it}
\thanks{}
\date{\today}
\maketitle

\begin{abstract}
By using the notion of contraction of Lie groups,
 we transfer $L^p-L^2$
estimates for  joint spectral projectors
from 
the unit complex sphere 
$\sfera$ in $\CC^{n+1}$   to the reduced Heisenberg group $h^{n}$.
In particular, we
deduce some estimates recently obtained by H. Koch and F. Ricci on $h^n$.
As a consequence,  we prove, 
in the
 spirit of Sogge's work, 
 a discrete restriction theorem for the sub-Laplacian $L$ on $h^n$.
\end{abstract}

\section{Introduction}
In the last twenty-five  years the notion of {\it{contraction}} (or {\it{continuous deformation}})  of Lie algebras and Lie groups,
introduced  in 1953 in a  physical context by E. In\"onu and E. P. Wigner,
was developed in a mathematical framework as well.
The basic idea is that,
given a Lie algebra $\frak{g}_1$,
from a family of non-degenerate transformations of
its structure constants it is possible to obtain,
in a limit sense,
a non-isomorphic 
 Lie algebra $\frak{g}_2$.

It turns out 
 that the deformed algebra $\frak{g}_2$ inherits
 analytic and geometric properties from
 $\frak{g}_1$
 and that the same holds for the corresponding Lie groups.
As a consequence, transference results have attracted considerable attention, in particular in the context of
Fourier  multipliers. In fact,  contraction has been successfully used
to transfer $L^p$  multiplier theorems from one Lie group to another one.
There is an extensive literature on such topic, centered about
deLeeuw's theorems;
we only mention here
the results by 
 A. H. Dooley, G. Gaudry, J. W. Rice and  R. L. Rubin
 (   \cite{Dooley}, \cite{DooleyGaudry},
 \cite{DooleyRice1}, \cite{DooleyRice2}, \cite{Rubin}),  concerning, in particular,  
contraction of  rotation groups and semisimple Lie groups.
  \vskip0.2cm
  The primary purpose of this paper  is
to show that contraction
is an effective tool to transfer  
$L^p$ eigenfunction bounds as well. In particular, we shall focus on  a contraction 
 from the complex unit sphere
$\sfera$ in $\CC^{n+1}$   to the reduced Heisenberg group $h^{n}$.

\vskip0.2cm
We recall that, if $P$ is a second order 
self-adjoint 
elliptic 
differential operator on a compact manifold $M$
and if $P_{\lambda}$ denotes the spectral
projection corresponding to the eigenvalue $\lambda^2$, 
a classical problem is to estimate the norm $\nu_p$ 
of $P_{\lambda}$  as an operator from $L^p (M)$, $1\le p\le 2$,  to $L^2(M)$.
Sharp estimates for $\nu_p$ have been obtained by C. Sogge (\cite{Sogge93}), who proved that
\begin{equation}\label{sogge}
||P_{\lambda}||_{(p,2)}
\le
C\lambda^{\gamma(\unpi,n)}\,\; 1\le p\le 2\,,
\end{equation}
where $\gamma$
is the piecewise
affine function on $[\unme,1]$
defined by
$$\gamma(\unpi,n):=
\begin{cases}\!
{n\left(\unpi-\frac{1}{2}\right)-\frac{1}{2} 
}
&\,\text{if
  $1\le p\le\tp$}\cr
{\frac{n-1}{2}(\frac{1}{p}-\unme)}
&\,\text{if
  $ \tp\le p\le 2$,}\cr
\end{cases}
$$
with  {\it{critical point}}
$\tp$  given by $\tp:=
  { 2\frac{n+1}{n+3}}$.
  
The starting point for our approach is a sharp two-parameter estimate for joint
spectral projections on complex spheres, recently obtained by the first author (\cite{Casarino2}). More precisely, 
we consider
the Laplace-Beltrami operator
$\Delta_{S^{2n+1}}$ and the Kohn Laplacian $\cL$ on
$S^{2n+1}$
(this set yields a basis for
the algebra of $U(n+1)$-invariant differential operators on
$S^{2n+1}$).
It is  possible  to work out 
a joint spectral theory. In particular, we denote by
$\hll$, $\ell,\ell'\ge 0$,  the joint  eigenspace  with eigenvalue
$\mu_{\ell,\ell'}$ for $\Delta_{S^{2n-1}}$, where
$\mu_{\ell,\ell'}:=
-\somma\left( \ell+\ell'+2n-2\right)$, and  with  eigenvalue
$\lambda_{\ell,\ell'} $ for $\cL$, where $\lambda_{\ell,\ell'}:=
-2\ell\ell'-(n-1)\somma$ (\cite{Klima}).
It is a classical fact 
([VK, Ch.11])
that
\begin{equation}\label{decomposition}
L^2 \left(
S^{2n+1}\right)=\dis\sum_{\ell,\ell'=0}^{+\infty}\!\oplus\cH^{\ell\ell'}\,,
\end{equation}
where the series on the right  converges in the $L^2$-norm.

\noindent By the symbol $\pi_{\ell\ell'}$
we denote the
joint spectral projector from
$ L^{2}(S^{2n-1})$
onto
$\cH^{\ell\ell'}$.
In \cite{Casarino2}
we proved the following two-parameter
$L^p$ eigenfunction bounds
 \begin{equation}\label{alfabeta}
 ||\pi_{\ell,\ell'}||_{(p,2)}\lesssim
C \,\left( 2\ql+n-1\right)^{\alpha(\unpi,n) }
\left(1+\Ql\right)^{\beta(\unpi,n) }\,
\,\text{ for all $\ell,\ell'\ge 0$\,,}
\end{equation}
 where
 $\Ql:=\max \{\ell,\ell'\}$, $\ql:=\min \{\ell,\ell'\}$
 and
 $\alpha$ and $\beta$ are the piecewise affine functions
represented in Figure 1 at the end of Section 2.
We remark that the critical exponent is in our case  $\frac{2(2n+1)}{2n+3}$
and 
cannot be directly deduced from  Sogge's results. 
Observe moreover that $2\ql+n-1$ and $\Ql$ are related to the
 eigenvalues
 $\lambda_{\ell,\ell'}$ and
 $\mu_{\ell,\ell'}$, since
 they  grow, respectively, as
$\frac{\left|\lambda_{\ell,\ell'}\right|}{\ell+\ell'}$ and
$\left|\mu_{\ell,\ell'}\right|^{\unme }$.
\vskip0.2cm

On the other hand, 
on the  reduced Heisenberg group
$\heisenunouno$,  defined as
$\heisenunouno:=
\CC
^n\times\TT$, with product
$$(\underz,e^{it})(\underw,e^{it'}):=
\left(
\underz+\underw,e^{i\left(t+t'+ \Im m\, \underz\bar{\underw}\right)}
\right)\,,$$
with $\underz,\underw\in\CC
^{n}$,
$t,s\in\RR
$,
we consider
the sub-Laplacian $L$ 
and the operator 
$i^{-1}\partial_t$.
The pairs $(2|m|(2k+1), m)$, with $m\in\ZZ\setminus\{0\}$
and $k\in\NN$, give the discrete joint spectrum of these operators.
Recently H. Koch and Ricci  proved the following 
$L^p-L^2$ estimate
for the 
orthogonal projector
$P_{m,k}$ onto the joint eigenspace 
\begin{equation}\label{Kochriccistima}
 ||P_{m,k}||_{\left(L^{p}(h^{n}), L^2 (h^{n})
\right) }\lesssim
C \,{\left(2k+n\right)}^{\alpha(\unpi,n) }\cdot
|m|^{\beta(\unpi,n)}\,,
\end{equation}
$1\le p\le2$, 
where
$\alpha$ and $\beta$ are given by
(\ref{alfabeta}) (\cite{KochRicci}).

We start  showing in Section 2 that
$P_{m,k}$ may be obtained as limit in the $L^2$-norm
of a sequence of  joint spectral projectors on $\sfera$.
Then  we give an alternative proof
of  $(\ref{Kochriccistima})$
 by a contraction argument.

A contraction from $SU(2)$ to the one-dimensional Heisenberg group $H^1$ was studied by F. Ricci and Rubin (\cite{Ricci}, \cite{RicciRubin}).
In \cite{Casarino2} the first author used some ideas from \cite{Ricci} 
to transfer $L^p-L^2$ estimates for norms of harmonic projection operators
from the unit sphere $S^3$
in $\CC^2$ to the reduced Heisenberg group 
$h^1$.
In this paper we discuss the higher-dimensional case.

A contraction from the unit sphere $\sfera$  to the Heisenberg group $H^n$
for $n>1$  was analyzed by Dooley and S. K. Gupta;
in a first paper 
they  adapted the notion of Lie group contraction to  
the  homogeneous space
$U(n+1)\// U(n)$  and described  the relationship between 
certain unitary  irreducible representations of $U(n+1)$ and $H^n$
(\cite{DooleyGupta1}),  in a second paper  
they  proved
a deLeeuw's type theorem on $H^n$ by transferring results  from $\sfera$
(\cite{DooleyGupta2}).
The contraction we use here is essentially that introduced by Dooley and Gupta;
anyway, their approach is mainly  algebraic, while our interest is adressed 
to  the analytic features of the problem.  

\vskip0.2cm
As an application of 
$(\ref{alfabeta}$)
 we prove in Section 3 a discrete restriction theorem for the sub-Laplacian $L$ on $h^n$ in the
 spirit of Sogge's work (\cite{Sogge}, see also (\ref{sogge})).
 More precisely, let
$Q_N$ be the spectral projection  corresponding to the eigenvalue $N$ associated to $L$ on $h^n$,
that is
$$Q_N f:=\displaystyle\sum_{(2k+n)|m|=N}P_{m,k}f\,.$$
 
\noindent
The study of $L^p-L^2$ mapping properties of $Q_N$  was suggested by D. M\"uller in his paper
about the restriction theorem on the Heisenberg group  (\cite{Mueller}).
In  \cite{Thangavelu1}  Thangavelu proved that
\begin{equation}\label{thangavelu}
\vert\vert Q_N\vert\vert _{\left(
L^{p}(h^{n}), L^2 (h^{n})\right)}
 \le C \,\left(N^n d(N) \right)^{\unpi-\unme} \,,\qquad
 1\le p\le 2\,,
\end{equation}
where $d(N)$ is the divisor-type function
defined  by
\begin{equation}
\label{defdN}
d(N):=\displaystyle
\sum_{2k+n | N} \frac{1}{2k+n}\,,
\end{equation}
and the estimate is sharp 
for $p=1$.
By  $a|b$ we mean that $a$ divides $b$.
Other types of restriction theorems on the Heisemberg group were
discussed by
Thangavelu in \cite{ThangaveluStudia}.

By using orthogonality,
 we add up the estimates
in 
$(\ref{alfabeta}$)
and obtain $L^p-L^2$ bounds for the norm of 
 $Q_N$, which 
in some cases improve
 (\ref{thangavelu}).
The exponent appearing in (\ref{thangavelu})
is an affine function of $\unpi$.
In our estimate the exponent of $d(N)$ is, like in Sogge's results,
a piecewise affine function of $\unpi$.
In other words,  there is
a critical point $\tp$ where the slope 
of the exponent changes. 
This critical point is the same that was found 
on complex spheres 
 (\cite{Casarino2}).

Our bounds are in general not sharp.
The reason  is that with our procedure
we disregard the interferences
between eigenfunctions.
We show however
that there are arithmetic progressions
${N_{m}}$ in $\NN$ for which our estimates 
 for 
 $\vert\vert Q_{N_{m}}\vert\vert _{\left(
{p}, 2 \right)}$
are sharp and better than 
(\ref{thangavelu}).
Moreover, since
 the behaviour of $d(N)$ is highly irregular, 
we inquire about the average size of 
$\vert\vert Q_{N}\vert\vert _{\left(
{p}, 2 \right)}$.
We prove in this case that $L^p-L^2$ estimates 
do not involve 
divisor-type functions and that the
critical point disappears.

It is a pleasure to thank Professor Fulvio Ricci for his valuable help.

\section{Preliminaries}

In this section we introduce some notation and  recall a few results, that will   be used 
 in the following.

\subsection  {\it{Some notation}}
For $n\ge 1$ let ${\CC}^{n+1}$ denote the n-dimensional complex space
endowed with the scalar product
$<\underz,\underw>:=\zu \wub +\ldots+\zn\wnb$,  $\underz,\underw \in{\CC}^{n+1}$, and
let $S^{2n+1}$ denote  the unit sphere in
$\CC
^{n+1}$, that is $$ S^{2n+1}
:=\{\underz=(\zu,\ldots,\zn)\in{\CC
}^{n+1}\,:\,
<\underz,\underz>=1\}\,.$$
The symbol $\uno$
will denote the north pole of
$S^{2n+1}$, that is
$\uno:=(0,\ldots,0,1)$.
\par
\vskip0.2cm

For every $\ell,\ell'\in\NN$
the symbol $\hllp$
will denote the space of the restrictions to
$S^{2n+1}$  of harmonic polynomials
 $p(\underz,\zetab)=
p(\zu,\ldots,\zn, \zub,\ldots,\znb)$,
of homogeneity degree $\ell$ in
$\zu,\ldots,\zn$ and of homogeneity degree $\ell'$ in
$(\zub,\ldots,\znb)$, $i.e.$ such that
$$p(a\underz,b\zetab)=
a^{\ell}b^{\ell'}p(\underz,\zetab)\,,\;\;a,b\in\RR\,,\;\;	\underz\in\CC^n
\,.$$
For a detailed description of the spaces $\hllp$ see Chapter 11 in \cite{KlimykVilenkin}.
We only recall here that 
a polynomial $p$ in $\underz,\zetab$ is said to be harmonic if
\begin{equation}
\label{defLaplace}
\Delta_{\sfera} p:= \frac{1}{4}\big(
\frac{\partial^{2}}{\partial\zu\partial\zub}
+ \ldots + \frac{\partial^{2}}{\partial\zn\partial\znb}
\big)p=0,
\end{equation}
where $\Delta_{\sfera}$ denotes the Laplace-Beltrami operator.
\vskip0.2cm
A zonal function of bidegree $(\ell,\ell')$
on $S^{2n+1}$ is a function 
in $\hllp$, which
is constant on the orbits of the
stabilizer of $\uno$ (which is isomorphic to $U(n)$).
Given a zonal function $f$, 
we may associate to  $f$
a map ${}^b\!{f}$  on the unit disk by 
$$f(\underz)={}^b\!{f}(<\underz,\uno>)\;,\,\underz\in
\sfera\,,$$
(by using the notation in  Section 11.1.5 of \cite{KlimykVilenkin} 
we have 
$<\underz,\uno>=z_{n}=e^{i\varphi}\cos\theta$,
where
$\varphi\in [0,2\pi]$ and $\theta\in [0,\frac{\pi}{2}]$).
\par
By means of ${}^b\!{f}$ we may define a 
convolution between a
zonal function $f$  and an arbitrary function
$g$ on $\sfera$. More precisely, we set
$$\big(f*g\big)(\underz):=\int_{\sfera}{}^b{f}(<\underz,\underw>)g(\underw)
d\sigma(\underw)\,,$$
where $d\sigma$
is the measure invariant under the action of the unitary group $U(n+1)$
(see ($\ref{misurasfera}$) 
for an explicit formula).
In the following we shall write $f(\theta,\varphi)$ instead of
${}^b{f}(e^{i\varphi}\cos\theta)$.
\vskip0.3cm
Let  $L^{2}(S^{2n+1})$ be
the Hilbert space of functions on
$S^{2n+1}$
endowed with the inner product
$(f,g):= \int_{S^{2n+1}}
f(\underz){\overline{g(\underz)}}d\sigma (\underz)$.
\par
It is a classical fact (\cite{KlimykVilenkin},
Ch. 11)
that $L^2(\sfera)$ is the direct sum of the pairwise orthogonal 
and $U(n+1)$-invariant
subspaces $\hllp$, $\ell,\ell'\ge 0$.
In other words,
every $f\in L^2 (\sfera)$
admits a 
unique expansion
$$f=\dis\sum_{\ell,\ell'=0}^{+\infty}\!
Y^{\ell\ell'}\,, 
$$
where $Y^{\ell\ell'}
\in\hllp$ for every $\ell,\ell'\ge 0$ and
the series at the right converges to $f$ in
the  $L^{2}(S^{2n+1})$-norm.
\par
The orthogonal  projector onto $\hllp$
\begin{equation}\label{definizioneproiettore}
\pllp\,: L^{2}(S^{2n-1})\ni f\mapsto
Y^{\ell\ell'}  \in\hllp\,\end{equation}
may be written as 
$$\pllp f:={}^b\zll *f\,,$$
where $\zll$ is the zonal function from $\hllp$, given by
\begin{equation}\label{zonali}
\begin{aligned}
^{b}\!\zll (\theta,\varphi)&:=
\frac{\dll}{\omega_{2n+1}}
\frac{\ql! (n-1)!}{(\ql+n-1)!}
e^{i(\ell'-\ell)\varphi}
    (\cos\theta)^{|\ell-\ell'|}
P_{\ql}^{(n-1,|\ell-\ell'|)}(\cos 2\theta)\;\;\\
&\qquad\qquad
\,\ell\,,\ell'\ge 1,
\,\varphi\in [0,2\pi]\,, \,\,\theta\in [0,\frac{\pi}{2}]\,.
\end{aligned}
\end{equation}
where
$\ql=\min \big( \ell,\ell'\big)$,
${\omega_{2n+1}}$ denotes the surface area of
${S^{2n+1}}$,
$P_{\ql}^{(n-1,|\ell-\ell'|)}$ is the Jacobi polynomial and
$$\dll:=\text{dim}\hll=
n\cdot \frac{\ell+\ell'+n}{\ell\ell'}
\begin{pmatrix} \ell+n-1\\
	 \ell-1\end{pmatrix}
\begin{pmatrix} \ell'+n-1\\
	 \ell'-1\end{pmatrix}\,\text{ for all $\ell,\ell'\ge 1$.}$$

Recall finally that 
 $\hlo$ consists of holomorphic polynomials
 and
$\cH^{0,\ell}$ consists of polynomials whose complex conjugates are
holomorphic.
In both cases, the dimension of the space is given by
$${\text{dim}\,\hlo}=
{\text{dim}\,\cH^{0,\ell}}=
{
{ \binom{\ell+n-1}{\ell}}}$$
and the zonal function is
$$\zlo (\theta,\varphi):=\frac{1}{\omega_{2n-1}}
\binom{\ell+n-1}{\ell}
e^{-i\ell\varphi}
    (\cos\theta)^{\ell}
,\;\;
\varphi\in [0,2\pi],\,\,\theta\in [0,\frac{\pi}{2}].$$
\vskip0.2cm
In this paper we shall adopt the convention that $C$ denotes a constant which is not necessarily the same at each occurrence.
\vskip0.2cm

\subsection{}{\it{Some useful results.}}
In order to transfer  $L^p$ bounds from $\sfera$ to $\heisenunouno$
we shall need
both a pointwise estimate for the Jacobi polynomials, due to Darboux and Szeg\"o
([Sz, pgs. 169,198]), 
and a
Mehler-Heine-type formula, relating Jacobi and Laguerre polynomials
(\cite{Szego}, \cite{Ricci}). 

\begin{lemma}\label{stimepuntualiszego }
Let $\alpha,\beta>-1$.
Fix $0<c<\pi$.
Then
$$
 P^{(\alpha,\beta)}_{\ell}(\cos\theta)=
 \begin{cases}\!
O\left(
{\ell^{ \alpha}}\right)
&\,\text{if
  $0\le \theta\le {\frac{c}{\ell}} $, }\cr
\ell^{-\unme}k(\theta)
\left(
\cos\left(N_{\ell}\theta+\gamma\right)
+\left(\ell\sin\theta\right)^{-1}O(1)\right)
&\,\text{if $
{\frac{c}{\ell}}\le\theta\le\pi-{\frac{c}{\ell}}$}\cr
O\left(
{\ell^{ \beta}}\right)
&\,\text{if
  $\pi-{\frac{c}{\ell}}\le\theta\le\pi$, }\cr
\end{cases}$$
where $k(\theta):=
\pi^{\unme}
\left( \sin {\frac{\theta}{2}}
\right)^{-\alpha-\unme}
\left( \cos {\frac{\theta}{2}}
\right)^{-\beta-\unme}$,
$N_{\ell}:=
\ell+
 {\frac{\alpha+\beta+1}{2}}$,
 $\gamma:=
 -(\alpha+\unme)\pimezzi$.
 \end{lemma}
 
 \begin{proposition}\label{Mehler}
 {\rm{[R, pg.224]}}
Let $n\ge 1$ and let $x$ be a real number.
Fix $k$ and $j$ in $\NN$, $j\ge k$.  Then
\begin{equation}
\begin{aligned}
\lim_{N\to +\infty}
&\,\cos^{N-j-k} \left(\frac{x}{\sqrt{N-j-k}}\right)\cdot
P^{(j-k,N-j-k)}_{k}\left( \cos \frac{2x}{\sqrt{N-j-k}}\right)
\\
&
\qquad\qquad=\,
L^{j-k}_k \left(  x^2\right)\cdot
e^{-\frac{1}{2}\,x^2}\,.
\end{aligned}
\end{equation}
\end{proposition}

Our proof is based on the following two-parameter estimate
for the $L^p-L^2$ norm of the complex harmonic projectors
$\pi_{\ell,\ell'}$, defined by
(\ref{definizioneproiettore}).
\begin{theorem}
\label{teoremariassuntivo}
 {\rm{\cite{Casarino2}}}
Let $n\ge 2$ and let $\ell,\ell'$ be non-negative integers.
Then
 \begin{equation}\label{stimaduep}
||\pi_{\ell,\ell'}||_{(p,2)}\lesssim
C \,\left( \frac{2\ell\ell'+n(\ell+\ell')}{\ell+\ell'}\right)^{\alpha(\unpi,n) }
\left(\ell+\ell' \right)^{\beta(\frac{1}{p},n)}
\,\text{if
  $1\le p\le 2$,}
\end{equation}
where 
 \begin{equation}\label{alfadef}
\alpha(\unpi,n):=
\begin{cases}\!
{n\left(\unpi-\frac{1}{2}\right)-\frac{1}{2} 
}
&\,\text{if
  $1\le p< \tp$}\cr
{\frac{1}{4}-\frac{1}{2p}}
&\,\text{if
  $\tp\le p\le 2$,}\cr
\end{cases}\end{equation}
with $\tp=
2
{ \frac{2n+1}{2n+3}}$, 
and
 \begin{equation}\label{betadef}
\beta(\unpi,n)=
n\left( \unpi-\unme\right)
\,\text{for all 
  $1\le p\le 2$, }
 \end{equation}
The above estimates are sharp.

\end{theorem}

 \vspace*{2.4in}
$$
\begin{picture}(20,20)
\setcoordinatesystem  units <.05in,.05in> point at 0 0
\setplotarea x from -360 to 360, y from -360 to 360
\put(-55,15){\vector(0,1){180}}     

 \put(-55,40){\vector(1,0){130}}     

\put(70,25) {$\frac{1}{p}$}

\put(-10,37.77){\line(0,1){5}}    
\put(-19.12,10.77){$\frac{2n+3}{2(2n+1)}$}    
\put(15.12,36.77){\line(0,1){6}}    


\put(-58,33){\line(1,0){6}}    
\put(-102.12,22.77){$-\frac{1}{2(2n+1)}$}    

\put(-58,169){\line(1,0){6}}    
\put(-80.12,167){$\frac{n}{2}$}    

\put(-58,125.5){\line(1,0){6}}    
\put(-80.12,123){$\frac{n-1}{2}$}    

\put(-18.5,36.77){\line(0,1){6}}    

\thicklines
\put(70,180){\line(1,0){22}}
\put(93,178){$\beta$}
\thinlines
\put(70,168){\line(1,0){22}}    
\put(93,166){$\alpha$}



\thicklines
\put(-17.55,40.1){\line(1,4){32}}    
\thinlines
\put(-9.1,34.5){\line(1,4){23}}         
\put(-9.1,34.5){\line(-3,2){10}}         


\put(-73,40) {$O$}
 \put(12,25) {$1$}
  \put(-25,25) {${\frac{1}{2}}$}

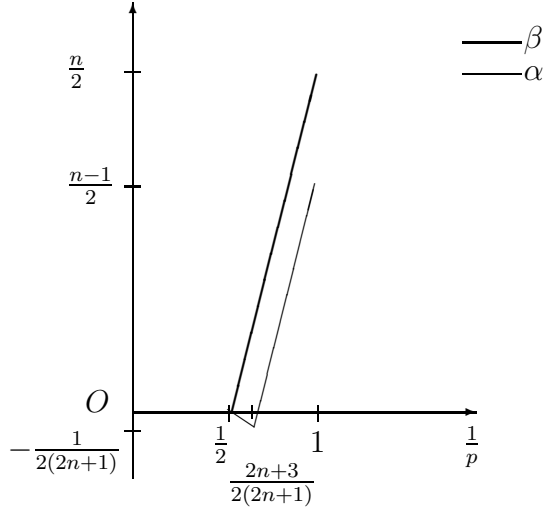
\put(-150,-40){Figure 1. The exponents $\alpha$ and $\beta$ as functions
of $\frac{1}{p}$}
\end{picture}
$$
\vskip2cm

\section{$L^p$ eigenfunction bounds on $\heisenuno$}
The Heisenberg group $H^n$ is a Lie group with underlying manifold $\CC
^n\times \RR
$,
endowed with the product
$$(\underz,t)(\underw,s):=
\left(
\underz+\underw,t+s+ \Im m\, \underz \cdot \overline{\underw}
\right)\,,$$
with $\underz, \underw\in\CC
^{n}$,
$t,s\in\RR
$.

\noindent We denote an element  in $H^1$
by
$\left(\rho{e^{i\varphi}},t\right)$, where
$\rho\in [0,+\infty)$, $\varphi\in[0,2\pi]$, $t\in\RR
$, 
and
an element  in $H^{n}$
by
$\left(\rho{\underline{\bf{\eta}}},t\right)$, where
$\rho\in [0,+\infty)$, $t\in\RR
$ and $\eeta\in S^{2n-1}$ is
given by
\begin{equation}
\label{coordeeta}
\eeta=
\begin{cases} {
e^{i\varphi_1}
\sin\theta_{n-1}
\sin\theta_{n-2}
\ldots
\sin\theta_{1}}\\
e^{i\varphi_2}
\sin\theta_{n-1}
\sin\theta_{n-2}
\ldots
\cos\theta_{1}\\
\vdots\\
e^{i\varphi_{n}}
\cos\theta_{n-1}\,,
\end{cases}
\end{equation}  with
$\varphi_k\in[0,2\pi]$, $k=1,\ldots,n$,
 and $\thetaj\in [0,\pimezzi]$, $j=1,\ldots, n-1$.
 
 \noindent
Observe that
$\eeta=\eeta\left({\Theta_{n-1}},\Phi_{n}\right)$, where
${\Theta_{n-1}}:=\left(
{\theta_{1}}\,,
{\theta_{2}}\,,
\ldots,
\theta_{n-1}\right)$
and
$\Phi_{n}:=\left(
\varphi_{1},
\ldots,
\varphi_{n}\right)$.

\noindent Define now a map
$\Psi:\heisenuno\to\sfera$
by
\begin{equation}
 \Psi\,:
 \, \left(
 \rho{\eeta},t\right)
\mapsto
(
{\Theta_{n-1}}\,,
\rho\,,\Phi_{n}
,t)\,,
\end{equation}
where
$
(
{\Theta_{n-1}}\,,
\rho\,,\Phi_{n},
t)\in
\sfera$ is given by
\begin{equation}
\label{coordsfera}
(
{\Theta_{n-1}}\,,
\rho\,,
\Phi_{n},
t):=
\begin{cases} {
e^{i\varphi_1}
\sin \rho\, \sin\theta_{n-1}
\, \sin\theta_{n-2}
\ldots
\sin\theta_{1}}\\
e^{i\varphi_2}
\sin \rho\,\sin\theta_{n-1}
\sin\theta_{n-2}
\ldots
\cos\theta_{1}\\
\vdots\\
e^{i\varphi_{n}}
\sin \rho\, \cos\theta_{n-1}\,\\
e^{it}\cos\rho\,.
\end{cases}
\end{equation}

\noindent We introduce in this way a coordinate system 
$
(
{\Theta_{n-1}}\,,
\rho\,, \Phi_{n},
t)$
on
$\sfera$, if $\rho$ and $t$ are
restricted, respectively, to $[0,\pimezzi]$ and $[-\pi,\pi]$.

\vskip0.4cm
The invariant measure $d\sigma_{\sfera}$ on $\sfera$ in the spherical coordinates
($\ref{coordsfera}$)
is
\begin{equation}\label{misurasfera}
\frac{n!}{2\pi^{n+1}}
\displaystyle
\Pi_{k=1}^{n}
d\varphi_{k}\,
\,dt\,\,
\sin^{2n-1}\rho\cos\rho\,
d\rho\,\,
\displaystyle
\Pi_{j=1}^{n-1}
\sin^{2j-1}\theta_j
\cos\thetaj\,d\thetaj.
\end{equation}
\noindent The factor 
$\frac{n!}{2\pi^{n+1}}
$
is introduced in order to make  the measure of the whole sphere equal to
$1$.


\noindent The Haar measure on $H^n$  in these coordinates is
$$
\frac{n!}{2\pi^{n+1}\sqrt{\omega_{2n+1}}}
\rho^{2n-1}
d\rho\,
d\varphi_{1}\ldots
d\varphi_{n}\,\displaystyle
\Pi_{j=1}^{n-1}
\sin^{2j-1}\theta_j
\cos\thetaj\,d\thetaj.$$
\vskip0.2cm

\vskip0.4cm
The {\it{reduced Heisenberg group}}
$\heisenunouno$ is defined as
$\heisenunouno:=
\CC
^n\times\TT$, with product
$$(\underz,e^{it})(\underw,e^{it'}):=
\left(
\underz+\underw,e^{i\left(t+t'+ \Im m\, \underz\bar{\underw}\right)}
\right)\,,$$
with $\underz,\underw\in\CC
^{n}$,
$t,s\in\RR
$.

\noindent Let now $f$ be a function on $\heisenunouno$, with compact support.
Let $\tif$ be the function $f$ extended by periodicity on $\RR
$ with respect  to
the variable $t$.
Define the function $\fnu$ on $\sfera$
by \begin{equation}\label{fnuvale}
f_{\nu}
(
{\rho}\,,
{\Theta_{n-1}}\,,
\Phi_{n},t)
:=
\nu^{n}\, \tif
(
\rho\sqrt{\nu}\,
\eeta,t\nu\,)
\,,\,\,\,\nu\in\NN.
\end{equation}

\begin{lemma}
\label{RicciRubin2}
Let $f$ be an integrable function on $\heisenunouno$ with compact support.
If $1\le p\le +\infty$, then
$$\nu^{-\frac{n}{p'}}
||f_{\nu}||_{L^p (\sfera)}<
||f||_{L^p (\heisenunouno)}\,\,\,\, \text{and}$$
$$
\displaystyle\lim_{\nu\to +\infty}
\nu^{-\frac{n}{p'}}||f_{\nu}||_{L^p (\sfera)}=
||f||_{L^p (\heisenunouno)}\,.$$
\end{lemma}
\begin{proof}
The proof is similar to that of Lemma 2 in \cite{RicciRubin}
and is omitted. Compare also with Lemma 4.3 in \cite{DooleyGupta2}.
\end{proof}

Throughout the paper we shall consider  a pair of strongly commuting operators on $h^n$.
The first  is 
the left-invariant sub-Laplacian $L$, defined by 
$$L:=
-
\displaystyle\sum_{j:=1}^n
\left(
X_j^2+
Y_j^2\right)\,,$$
where
$X_j:=
\partial_{x_{j}}
-y_{j}\partial_{t}\,$
and 
$Y_j:=
\partial_{y_{j}}
+x_{j}\partial_{t}\,$. 
The second  is the operator
$T:=i^{-1}\partial_t$.
These operators generate the algebra of  differential operators on $h^n$
invariant under  left translation
and under the action of the unitary group. 
One can  work out a joint spectral theory; 
the pairs $(2|m|(2k+n), m)$, with $m\in\ZZ\setminus\{0\}$
and $k\in\NN$,  give the discrete joint spectrum of 
$L$ and
$i^{-1}\partial_t$.
We shall denote by
$P_{m,k}$ 
the orthogonal projector
onto the joint eigenspace.

By considering the Fourier decomposition of functions
in $L^2 (h^n)$ with respect to the central variable, we obtain an orthogonal decomposition
of $L^2 (h^n)$ as
$$L^2 (h^n)=
\HH_0 \oplus
\HH\,,$$
where $\HH$ is given by
$$\HH:=
\{f\in L^2 (h^n)\,:\,
\int_{\TT} f(z,t)dt=0\,
\}\,.$$
The projectors $P_{m,k}$ map $L^2 (h^n)$ onto
$\HH$
and provide a spectral decomposition for $\HH$.
We point the attention on this decomposition, since the spectral analysis of $L$ on $\HH_0$
essentially reduces to the analysis of the Laplacian on $\CC^n$.

\vskip0.2cm
On the complex sphere $\sfera$ the algebra of $U(n+1)$-invariant differential operators 
is commutative and generated by two elements;
a basis is given by
the Laplace-Beltrami operator
$\Delta_{S^{2n+1}}$, 
defined by
($\ref{defLaplace}$),
and the Kohn Laplacian $\cL$ on
$S^{2n+1}$, defined by
$$\cL:=
\displaystyle
\sum_{j<k}M_{jk}\overline{M}_{jk}
+
{\overline{M}}_{jk}M_{jk}\,,$$
with
$$M_{jk}:=
\overline{z}_j
\partial_{z_k}
-
\overline{z}_k
\partial_{z_j}\,\;\;
\text{and}
\;\;\;
\overline{M}_{jk}:=
z_j
\partial_{\overline{z}_k}
-
z_k
\partial_{\overline{z}_j}
\,.$$
We shall call
$\hll$ the joint eigenspace of
$\Delta_{S^{2n+1}}$  and $\cL$, with eigenvalues
respectively
$\mu_{\ell,\ell'}:=
-\somma\left( \ell+\ell'+2n\right)$ and 
$\lambda_{\ell,\ell'}=
-2\ell\ell'-n\somma$ ([Kl]).
\vskip0.2cm
The next task is proving that  the joint spectral projection $P_{m,k}$ on $\heisenunouno$
 may be obtained as limit in the
$L^2$-norm
of an appropriate sequence of  joint spectral projectors on $S^{2n+1}$.

\begin{proposition}\label{lemmaapprox}
Let $f$ be a continuous  function on $\heisenunouno$, with compact support.
Take $m\in\NN\setminus\{0\}$ and $k\in\NN$.
For every $\nu\in\NN$ let $N(\nu)\in\NN$ be
such that
\begin{equation}
\label{condizioneFulvio1}
\lim_{\nu\to +\infty}
{\frac{N(\nu)}{{\nu} }}
=m\,.
\end{equation}
Then
\begin{equation}
\label{limitepkjn}
||P_{m,k}f||_{L^2 (\heisenunouno)}=\lim_{\nu\to +\infty}
{\frac{1}{{\nu}^{\nmenouno} }}
||\pkne\fnu||_{L^2 (\sfera)}
\,, \text{ and }
\end{equation}
\begin{equation}
\label{limitepjkn}
||P_{-m,k}f||_{L^2 (\heisenunouno)}=\lim_{\nu\to +\infty}
{\frac{1}{{\nu}^{\nmenouno} }}
||\pnke\fnu||_{L^2 (\sfera)}
\,.
\end{equation}
\end{proposition}

\begin{proof}
The scheme of the proof is similar to that
of Proposition 4.4 in 
\cite{Casarino2}.
Since  the higher dimensional case
is   more involved,
 we present 
 the proof
 for more transparency. 
 
Fix  two integers $m>0$ and  $k\in\NN$.

First of all, 
if $\underz,\underw\in\CC^n$,
by writing  $\underz:=\rho\eeta$ and
$\underw:=\rho'\eetap$,
with $\rho,\rho'\in[0,+\infty)$
and $\eeta$,$\eetap\in S^{2n-1}$,
a simple computation yields 
\begin{equation}\label{imma}
\begin{aligned}
\Im m (\underz\cdot\overline{\underw})=
\rho\rho'&\cdot
\left(
\sin (\varphi_1-\varphi'_1)
\sin\theta_{n-1}
\sin\theta'_{n-1}
\ldots\ldots
\sin\theta_1
\sin\theta'_1
\right.\\
&\left.
+\sin (\varphi_2-\varphi'_2)
\sin\theta_{n-1}
\sin\theta'_{n-1}
\ldots\ldots
\cos\theta_1
\cos\theta'_1
+
\ldots\right.\\
&\left.\ldots+
\sin (\varphi_{n}-\varphi'_{n})
\cos\theta_{n-1}
\cos\theta'_{n-1}\right)
\,\\
\end{aligned}
\end{equation}
and
\begin{equation}\label{modulo}
\begin{aligned}
|\underz-\underw|^2&=
\rho^2+\rho'^2
-2\rho\rho'
\cdot \left(
\cos (\varphi_1-\varphi'_1)
\sin\theta_{n-1}
\sin\theta'_{n-1}
\ldots
\sin\theta_1
\sin\theta'_1\right.\\
&\left.\qquad\qquad\qquad+
\cos (\varphi_2-\varphi'_2)
\sin\theta_{n-1}
\sin\theta'_{n-1}
\ldots
\cos\theta_1
\cos\theta'_1+\ldots\right.\\
&\left.\qquad\qquad\qquad\ldots+
\cos (\varphi_{n}-\varphi'_{n})
\cos\theta_{n-1}
\cos\theta'_{n-1}
\right)\,.
\\
\end{aligned}
\end{equation}

Now, by the symbol 
$\Phi_{k,k}^m$
we denote the joint eigenfunction
for $\cL$ and $i^{-1}\partial_t$
(for more details and an explicit expression see, for example,
[FH, Chapitre V]).
Orthogonality of joint spectral projectors yields
\begin{align*}
||&P_{m,k}f||_{L^2 (\heisenunouno)}^2=
<P_{m,k}f,f>_{L^2 (\heisenunouno)}
=\int_{\heisenunouno}
f*\Phi_{k,k}^m (\underz,t)\, \overline{f(\underz,t)}\,d\underz\,dt\\
&\qquad=\int_{\heisenunouno}
\!\int_{\heisenunouno}
\Phi_{k,k}^m
\left(
\underz-\underw, t-t'+
\Im m (\underz
\cdot\overline{\underw})\right)
f(\underw,t')
\,d\underw\,dt'\, \overline{f(\underz,t)}\,d\underz\,dt\\
&=
m^{n}
\int_{\heisenunouno}
\!\!\int_{\heisenunouno}
e^{ i\,m\left( t-t'+\Im m (\underz\cdot\overline{\underw})\right)}
L^{n-1}_k
\left(
m\, |\underz-\underw|^2
\right)
e^{
-\unme
m\, |\underz-\underw|^2
}
f(\underw,t')
\,d\underw\,dt'\\
&\qquad\qquad\qquad\qquad\qquad \overline{f(\underz,t)}\,d\underz\,dt\,.
\end{align*}

\par
Now we shall deal with the right-hand side in
($\ref{limitepkjn}$).
For the sake of brevity we set
$$d\Phi_{(n)}:=
d\varphi_1\,,\ldots,d\varphi_n\;\;\text{ and}$$
$$d\Theta_{(n-1)}:=
\Pi_{j=1}^{n-1}
\sin^{2j-1}\theta_j
\cos\thetaj\,d\thetaj.$$
From the orthogonality of the joint spectral projectors $\pi_{\ell,\ell'}$
in $L^2 (\sfera)$
and from ($\ref{fnuvale}$)
we deduce
\begin{align*}
||&\pkne \fnu||_{L^2 (\sfera)}^2=
<\pkne \fnu, \fnu>_{L^2 (\sfera)}\\
&=\int_{\sfera}
\left(\pkne \fnu\right)
 (
{   {\Theta_{n-1}}  }\,,
\rho\,,
\Phi_{n},t)
\,\overline{\fnu
(
  {\Theta_{n-1}}\,,
\rho\,,\Phi_{n}
,t)
}
\,
d\sigma_{\sfera}\\
&=
\frac{n!}{2\pi^{n+1}\,\nu}
\int_{A_{\nu}}
\left(
\pkne \fnu\right)
(
{\Theta_{n-1}}\,,
{\frac{\rho}{\sqrt{\nu}} }\,,
\Phi_{n}
,
{\frac{t}{{\nu}} })
\overline{
\tilde{f}
\left(
{\Theta_{n-1}}\,,
{\rho}\,,
\Phi_{n}
,t
\right)}
\left(\frac{\sin\,{{\frac{\rho}{\sqrt{\nu}} }}}{{\frac{\rho}{\sqrt{\nu}}
}}\right)^{2n-1}\!\!
\\
&\qquad\qquad
\cos{\frac{\rho}{\sqrt{\nu}} }\,
{\rho}^{2n-1}
d\rho\,
d\Theta_{(n-1)}
\,d\Phi_{(n)}\, dt\,\\
&=
\frac{n!^2}{4\pi^{2n+2}\nu^2}
\int_{A_{\nu}}
\bigg(
\int_{A_{\nu}}
{}^{b}\zkne
\left(
<(
{\Theta_{n-1}}\,,
{\frac{\rho}{\sqrt{\nu}} },
\Phi_{n}
,
{\frac{t}{{\nu}} }
)\,
,
(
{\Theta'_{n-1}}\,,
{\frac{\rho'}{\sqrt{\nu}} },
\Phi'_{n}
,
{\frac{t'}{{\nu}} }
 )>\right)
\\
&\qquad
\qquad
{\tilde{f}}
\left(
{\Theta'_{n-1}}\,,
{\rho'}
,\,
\Phi'_{n}
,
{t'}
\right)
\left(\frac{\sin\,{{\frac{\rho'}{\sqrt{\nu}} }}}{{\frac{\rho'}{\sqrt{\nu}}
}}\right)^{2n-1}\!\!
\cos{\frac{\rho'}{\sqrt{\nu}} }\,
{\rho'}^{2n-1}
d\rho'\,
d\Theta'_{(n-1)}\,d\Phi'_{(n)}\, dt'\bigg)\,\\
&\qquad\qquad\overline{
\tilde{f}
\left(
{\Theta_{n-1}}\,,
{\rho}
,\,
\Phi_{n}
,
{t}
\right)}
\left(\frac{\sin\,{{\frac{\rho}{\sqrt{\nu}} }}}{{\frac{\rho}{\sqrt{\nu}}
}}\right)^{2n-1}
\!\!\cos{\frac{\rho}{\sqrt{\nu}} }
{\rho}^{2n-1}
d\rho\,
d\Theta'_{(n-1)}\,d\Phi_{(n)}\, dt\,\\
\end{align*}
where
the integration set 
${A_{\nu}}$
is given by
\begin{equation}
\begin{aligned}
{A_{\nu}}&:=
\left\{
(
\rho\,,
{\Theta_{n-1}}\,,
\Phi_{n}
,t)\,
:\,
0\le\rho\le
{\frac{\pi}{2}\sqrt{\nu}}\,,
{0}\le\varphi_k\le {2\pi}\,,
k=1,\ldots,n\,,\,\right.\\
&
\left.\qquad \quad 0\le \thetaj\le
{\frac{\pi}{2}}\,,
j=1,\ldots,n-1\,,\,
{-\pi\nu}\le t\le {\pi\nu}
\right\}\,.
\end{aligned}
\end{equation}

Now by using 
$(\ref{coordsfera})$ we compute
the inner product in $\CC^{n+1}$
 \begin{align*}
<(
{\Theta_{n-1}}\,,
&{\frac{\rho}{\sqrt{\nu}}} ,
\Phi_{n-1},
{\frac{t}{{\nu}} }
)\,
,
(
{\Theta'_{n-1}}\,,
{\frac{\rho'}{\sqrt{\nu}} },
\Phi'_{n-1},
{\frac{t'}{{\nu}} }
 )>
=\\
&=
e^{i(\varphi_1-\varphi'_1)}
\sin(\rreps)
\sin(\rrepsp)
\sin{\theta_{n-2}}
\sin{\theta'_{n-2}}
\ldots
\sin{\theta_{1}}
\sin{\theta'_{1}}\\
&\qquad
+e^{i(\varphi_2-\varphi'_2)}
\sin(\rreps)
\sin(\rrepsp)
\sin{\theta_{n-2}}
\sin{\theta'_{n-2}}
\ldots
\cos{\theta_{1}}
\cos{\theta'_{1}}\\
&\qquad+\ldots
+e^{i(\varphi_{n-1}-\varphi'_{n-1})}
\sin(\rreps)
\sin(\rrepsp)
\cos{\theta_{n-2}}
\cos{\theta'_{n-2}}
\\
&\qquad+
e^{i(t-t')\repn}
\cos(\rreps)
\cos(\rrepsp)
\\
&=R_{\nu} e^{i\psi_{\nu}}\,,
\end{align*}
where
\begin{align*}
R_{\nu}&=
1-{\frac{1}{2\nu}}
\bigg(
\rho^2+\rho'^2
-2\rho\rho'
\left(
\cos (\varphi_1-\varphi'_1)
\sin\theta_{n-1}
\sin\theta'_{n-1}
\ldots
\sin\theta_1
\sin\theta'_1\right.\\
&\left.\qquad\qquad\qquad+
\cos (\varphi_2-\varphi'_2)
\sin\theta_{n-1}
\sin\theta'_{n-1}
\ldots
\cos\theta_1
\cos\theta'_1+\ldots\right.\\
&\left.\qquad\qquad\qquad\ldots+
\cos (\varphi_{n}-\varphi'_{n})
\cos\theta_{n-1}
\cos\theta'_{n-1}
\right)\bigg)
+o({\frac{1}{\nu}})
\,,\,\nu\to +\infty\,, \,\text{and}
\\
\end{align*}
\begin{align*}
\psi_{\nu}&=
\arctan\bigg(
{\frac{1}{\nu}}
\rho\rho'
\left(
\sin (\varphi_1-\varphi'_1)
\sin\theta_{n-1}
\sin\theta'_{n-1}
\ldots
\sin\theta_1
\sin\theta'_1
\right.\\
&\left.\qquad\qquad
+\sin (\varphi_2-\varphi'_2)
\sin\theta_{n-1}
\sin\theta'_{n-1}
\ldots
\cos\theta_1
\cos\theta'_1
+
\ldots\right.\\
&\left.\qquad\qquad\ldots+
\sin (\varphi_{n}-\varphi'_{n})
\cos\theta_{n-1}
\cos\theta'_{n-1}\right)+\frac{t-t'}{\nu}+o({\frac{1}{\nu}})\bigg)
\,\qquad\nu\to +\infty\,.\\
\end{align*}
Thus
as a consequence of
($\ref{imma}$) and ($\ref{modulo}$)
we have
$$
R_{\nu}=
\cos\left( {{\frac{1}{\sqrt\nu}}}|\underz-\underw|\right)+
o({{\frac{1}{\nu}}})\,\,\;\;\;\text{and }\quad
\psi_{\nu}=
{\frac{1}{\nu}}(t-t')
+{\frac{1}{\nu}}\Im m \underz\,\overline{\underw}+
o({\frac{1}{\nu}})\,,$$
so that  formula ($\ref{zonali}$)
for the zonal function yields
\begin{align*}
{}^{b}\zkne&
\left(
<(
{\Theta_{n-1}}\,,
{\frac{\rho}{\sqrt{\nu}} },
\Phi_{n}
,
{\frac{t}{{\nu}} }
)\,,(
{\Theta'_{n-1}}\,,
{\frac{\rho'}{\sqrt{\nu}} },
\Phi'_{n}
,
{\frac{t'}{{\nu}} }
 )>\right)\\&=
\frac{(N(\nu))^{n}}{\omega_{2n+1}}\,
e^{i\left(N(\nu)-2k\right){\frac{1}{\nu}}\left( t-t'+\Im
m\underz\,\bar{\underw}+o(1)\right) }
\left(  \cos \left(
{\frac{1}{\sqrt\nu}}|\underz-\underw|\right)\right)^{|N(\nu)-2k|}
\\
&\qquad
P^{(n-1,|N(\nu)-2k|)}_k \left(\cos \left(
{\frac{2}{\sqrt\nu}}|\underz-\underw|\right)\right)
+o({\frac{1}{\nu}})\,,\,\,\nu\to +\infty\,.
\end{align*}

By using condition  (\ref{condizioneFulvio1})
and the Mean Value Theorem,
we easily check that 
$$\frac{1}{\nu^{n}}
||\pkne \fnu||_{L^2 (\sfera)}^2
={\cI}_{\nu}^M+
\cI_{{\nu}}^R\,,$$
where the remainder term $\cI_{\nu}^R$ satisfies
$\displaystyle\lim_{\nu\to+\infty}\cI_{\nu}^R=0$, while the  main term $\cI_{\nu}^M$ is given by
\begin{align*}\label{principale}
&\cI_{\nu}^M
=\frac{n!^2}{4\omega_{2n+1}\pi^{2n+2}\nu^2}\,
\int_{A_{\nu}}
\left(
\int_{A_{\nu}}
\left({\frac{N(\nu)}{\nu}}\right)^{n}\,
e^{i
m\left( t-t'+\Im m\underz\,\bar{\underw}\right) }
\left(  \cos \left(
{\frac{1}{\sqrt\nu}}|\underz-\underw|\right)\right)^{|N(\nu)-2k|}
\right.\\
&\left.
\qquad\qquad
P^{(n-1,|N(\nu)-2k|)}_k
\left(\cos \left( {\frac{2}{\sqrt\nu}}|\underz-\underw|\right)\right)
{\tilde{f}}
\left(
{\rho'}
,\,
{\Theta'_{n-1}}\,,
\Phi'_{n}
,
{t'}
\right)
\left(\frac{\sin\,{{\frac{\rho'}{\sqrt{\nu}} }}}{{\frac{\rho'}{\sqrt{\nu}}
}}\right)^{2n-1}
\right.\,\\
&\left.\qquad\qquad
\cos{\frac{\rho'}{\sqrt{\nu}} }\,
{\rho'}^{2n-1}
d\rho'\,
d\Theta'_{(n-1)}\,
d\Phi'_{(n)}\,dt'\right)\,
\overline{
\tilde{f}
\left(
{\rho}
,\,
{\Theta_{n-1}}\,,
\Phi_{n}
,
{t}
\right)}
\left(\frac{\sin\,{{\frac{\rho}{\sqrt{\nu}} }}}{{\frac{\rho}{\sqrt{\nu}}
}}\right)^{2n-1}
\\
&\qquad\qquad
\cos{\frac{\rho}{\sqrt{\nu}} }\,
{\rho}^{2n-1}
d\rho\,
d\Theta_{(n-1)}\,d\Phi_{(n)}\, dt\,,\,\,\nu\to +\infty\,.
\\
\end{align*}
We shall now treat 
$\cI_{\nu}^M$
by means of  the Lebesgue Dominated Convergence Theorem.
First of all, we 
extend the integration set  in $\cI_{\nu}^M$,  
(this may be done, since
$f$ has compact support and  the integrand  is periodic
 with respect  to  $t$), 
 and we obtain 
 \begin{equation}
\label{termineprincipale1}
\begin{aligned}
\cI_{\nu}^M=&
\frac{n!^2}{4\pi^{2n+2}\omega_{2n+1}}\,
\int_{0}^{+\infty}\!
\int_{0}^{\pimezzi}
\ldots
\int_{0}^{\pimezzi}
\int_{0}^{2\pi}\!\!\!
\ldots
\int_{0}^{2\pi}
\int_{-\pi}^{\pi}\!\\
&\,\,\qquad\left(
\int_{0}^{+\infty}\!
\int_{0}^{\pimezzi}\!\!\!
\ldots\int_{0}^{\pimezzi}
\int_{0}^{2\pi}\!\!\!
\ldots
\int_{0}^{2\pi}
\int_{-\pi}^{\pi}\!
\left(
{\frac{N(\nu)}{\nu}}\right)^{n}
\quad e^{i\,m\left( t-t'-\Im m\underw\,\bar{\underz}\right) }
\,\right.\\
&\left.
\left(  \cos \left(
{\frac{1}{\sqrt\nu}}|\underz-\underw|\right)\right)^{|N(\nu)-2k|}
\,
P^{(n-1,|N(\nu)-2k|)}_k
\left(\cos \left( {\frac{2}{\sqrt\nu}}|\underz-\underw|\right)\right)
\,\right.\\
&\quad
\left.
{{f}}
\left(
{\rho'}
,\,
{\Theta'_{n-1}}\,,
\Phi'_{n}
,
{t'}
\right)
\left(\frac{\sin\,{{\frac{\rho'}{\sqrt{\nu}} }}}{{\frac{\rho'}{\sqrt{\nu}}
}}\right)^{2n-1}
\cos{\frac{\rho'}{\sqrt{\nu}} }\,
{\rho'}^{2n-1}
d\rho'\,
d\Theta'_{(n-1)}
d\Phi'_{(n)}\,
 dt'\right)
\\
&\quad
\overline{
{f}
\left(
{\rho}
,\,
{\Theta_{n-1}}\,,
\Phi_{n}
,
{t}
\right)}
\left(\frac{\sin\,{{\frac{\rho}{\sqrt{\nu}} }}}{{\frac{\rho}{\sqrt{\nu}}
}}\right)^{2n-1}
\cos{\frac{\rho}{\sqrt{\nu}} }\,
{\rho}^{2n-1}
d\rho\,
d\Theta_{(n-1)}
\,d\Phi_{(n)}\, dt\,.
\\
\end{aligned}
\end{equation}

By using Lemma
$\ref{stimepuntualiszego }$ and 
the Mehler-Heine formula as stated in
Lemma
$\ref{Mehler}$
(with $N=N(\nu)+j-k$, $j-k=n-1$ and
$x=\sqrt{\frac{N(\nu)-2k}{\nu}}|\underz-\underw|$),
 we may conclude as in Proposition 4.4 in \cite{Casarino2}.
 
The proof for ($\ref{limitepjkn}$) is completely analogous.
\end{proof}

\begin{theorem}
\label{teoremariassuntivo2}
Let $n>2$.
Take $m\in\ZZ\setminus\{0\}$ and $k\in\NN$.
Then
 \begin{equation}\label{stimaduep}
||P_{m,k}|| _{\left(L^{p}(h^{n}), L^2 (h^{n})
\right) }
\lesssim
\begin{cases}\!
C \,\left( 2k+n
\right)^{n\left(\frac{1}{p}-\frac{1}{2}\right)-\unme}
|m|^{n\left(\frac{1}{p}-\frac{1}{2}\right)}
&\,\text{if
  $1\le p< \tp$}\cr
C \,\left( 2k+n
\right)^{\frac{1}{4}-\frac{1}{2p}
}
 |m|^{n\left(\frac{1}{p}-\frac{1}{2}\right)}
&\,\text{if
  $\tp\le p\le 2$,}\cr
\end{cases}
\end{equation}
where $\tp=
2
{ \frac{2n+1}{2n+3}}$.
Moreover, the estimates are sharp.
\end{theorem}
\begin{proof}
Take $m>0$ (the other case being analogous).
For every $\nu\in\NN$ let $N(\nu)\in\NN$ be
such that
$$\lim_{\nu\to +\infty}
{\frac{1}{{\nu} }}
\cdot N(\nu)=m\,.
$$
Thus 
\begin{align*}||P_{m,k}f||_{L^2(\heisenunouno)}&=
\lim_{\nu\to +\infty}
{\frac{1}{{\nu^{\frac{n}{2}}}}}
||\pkne \fnu||_{L^2(\sfera)}\\
&
\le
\lim_{\nu\to +\infty}
\left({\frac{N(\nu)}{{\nu}}}\right)^{\frac{n}{2}}
{{{
\left( {\frac{2k\cdot (N(\nu)-k)}{N(\nu)}+n}\right)^{\frac{n}{2}}}}}
||\fnu||_{L^1(\sfera)}\\
&=
{{m}}^{\frac{n}{2}}
{{{(2k+n
)^{\frac{n-1}{2}}}}}
\lim_{\nu\to +\infty}|| \fnu||_{L^1(\sfera)}\\
&=
{{m}}^{\frac{n}{2}}
{{{(2k+n
)^{\frac{n-1}{2}}}}}
||f||_{L^1(\heisenunouno)}\,,
\end{align*}
where we used
first 
$(\ref{limitepkjn})$
and then
Theorem $\ref{teoremariassuntivo}$
and
Lemma
$\ref{RicciRubin2}$.

In the same way, we see that
\begin{align*}
||P_{m,k}f||_{L^2(\heisenunouno)}
&= \lim_{\nu\to +\infty}
{\frac{1}{{\nu^{\frac{n}{2}}}}}
||\pkne \fnu||_{L^2(\sfera)}\\
&
\le \lim_{\nu\to +\infty}
{\frac{1}{{\nu^{\frac{n}{2}}}}}
\left( {\frac{2k\cdot (N(\nu)-k)}{N(\nu)}+n}\right)^{-\frac{1}{2(2n+1)}}
(N(\nu))^{\frac{n}{2n+1}}
|| \fnu||_{L^{{2\frac{2n+1}{2n+3}}
}(\sfera)}
\\
&\le
(2k+n
)^{-\frac{1}{2(2n+1)}}
\lim_{\nu\to +\infty}
{\frac{1}{{\nu^{\frac{n}{2}}}}}
(N(\nu))^{\frac{n}{2n+1}}
\nu^{\frac{n(2n-1)}{2(2n+1)}}
|| f||_{L^{{2\frac{2n+1}{2n+3}}
}(\heisenunouno)}
\\
&
=
(2k+n
)^{-\frac{1}{2(2n+1)}}
m^{\frac{n}{2n+1}}
|| f||_{L^{{2\frac{2n+1}{2n+3}}
}(\heisenunouno)}
\,.\\
\end{align*}
An interpolation argument yields the thesis.
Finally, sharpness follows from arguments in \cite{KochRicci}.
\end{proof}

\section{A restriction theorem on $h^n$}

By applying  the bounds proved in
Section 2
we obtain
a restriction theorem 
for the spectral projectors associated to the sub-Laplacian $L$
on $h^n$.
Our theorem improves in some cases a previous result due to
 Thangavelu (\cite{Thangavelu1}).
 More precisely, let
$Q_N$ be the spectral projection  corresponding to the eigenvalue $N$ associated to $L$ on $h^n$,
that is
$$Q_N f:=\displaystyle\sum_{(2k+n)|m|=N}P_{m,k}f\,,$$
where $\pmk$ is the joint spectral projection operator introduced
in the previous section.
We  look for 
estimates of the type
\begin{equation}
\label{thangricci}
\vert\vert Q_N\vert\vert _{\left(
L^{p}(h^{n}), L^2 (h^{n})\right)}
\le C\,
N^{\sigma (p,n)} \,,
\end{equation}
for 
all $1\le p\le 2$, where the
exponent  $\sigma$ is in general a  convex function of  $\unpi$.
\vskip0.1cm
\noindent
In [Th91] Thangavelu proved that
\begin{equation}\label{thang}
\vert\vert Q_N\vert\vert _{\left(
L^{p}(h^{n}), L^2 (h^{n})\right)}
 \le C \,N^{n(\unpi-\unme)} d(N)^{\unpi-\unme}\,,\qquad
 1\le p\le 2\,,
\end{equation}
where $d(N)$ is the divisor-type function
defined  by
\begin{equation}
\label{defdN}
d(N):=\displaystyle
\sum_{2k+n | N} \frac{1}{2k+n}\,,
\end{equation}
and the estimate is sharp 
for $p=1$.
By  $a|b$ we mean that $a$ divides $b$.
\par
Thangavelu also proved that when $N=nR$,
with $R\in\NN$, 
then
$$
C\,N^{n(\unpi-\unme)} 
\le \vert\vert Q_N\vert\vert _{\left(
L^{p}(h^{n}), L^2 (h^{n})\right)}
\,,\qquad
 1\le p\le 2\,.
 $$
Here we show that there exist 
arithmetic progressions
${a_{N}}$ in $\NN$ such that
the estimate for $\vert\vert Q_{a_{N}}\vert\vert _{\left(
{p}, 2 \right)}
$ is sharp and better
 than
($\ref{thang}$) for $1<p<2$.

\begin{proposition}
\label{stimaQN}
Let $n\ge 1$.
Let $N$ be any positive integer number.

\noindent
Then
for every $1\le p\le 2$ 
 \begin{equation}\label{stimaQNp2}
||Q_N|| _{\left(L^{p}(h^{n}), L^2 (h^{n})
\right) }
\le
 C \,N^{n(\unpi-\unme)} d(N)^{\rho (\unpi,n)}\,,
\end{equation}
where $\rho$ is defined by 
 \begin{equation}\label{ro}
 \rho (\unpi,n):=
\begin{cases}\!
\unme
&\,\text{if
  $1\le p< \tp$}\cr
(2n+1)\left(  {\frac{1}{2p}-\frac{1}{4}
}\right)
&\,\text{if
  $\tp\le p\le 2$,}\cr
  \end{cases}
 \end{equation}
 with $\tp=
2
{ \frac{2n+1}{2n+3}}$,
and $d(N)$ is given by
($\ref{defdN}$).

\end{proposition}

\begin{proof}
For $p=1$ our estimate coincide with 
$(\ref{thang})$; nonetheless we give a different, simpler  proof:
\begin{align*}
||Q_N f||^2_{L^2(\heiserid)}&=
\displaystyle\vert\vert\sum_{(2k+n)|m|=N}
P_{m,k}f\vert\vert^2_{L^2(\heiserid)}=
\displaystyle\sum_{(2k+n)|m|=N}
\vert\vert
P_{m,k}f\vert\vert^2_{L^2(\heiserid)}\\
&\le C\,
\displaystyle\sum_{(2k+n)|m|=N}
{{m}}^{{n}}
{{{(2k+n)^{{n-1}}}}}
||f||^2_{L^1(\heiserid)}\,,
\\
&\le C
{{N^{n}}}
\displaystyle\sum_{2k+n|N}
{\frac{1}{2k+n}}||f||^2_{L^1(\heiserid)}\,,\\
\end{align*}
whence
\begin{equation}\label{stima1-2}
||Q_N ||_{\left( L^1, L^2\right) }\le CN^{\frac{n}{2}} \left( d(N)\right)^{\frac{1}{2}}  \,.
\end{equation}
For $p=2$ the bound is obvious, since $Q_N$ is an orthogonal projector.
Finally, for $p=\tp$ one has
\begin{align*}
||Q_N f||^2_{L^2(\heiserid)}&=
\displaystyle\sum_{(2k+n)|m|=N}
\vert\vert
P_{m,k}f\vert\vert^2_{L^2(\heiserid)}\\
&\le C
\displaystyle\sum_{(2k+n)|m|=N}
(2k+n)^{-\frac{1}{2n+1}}
|m|^{\frac{2n}{2n+1}}
||f||^2_{L^{\tp}(\heiserid)}\,,\\
&= C
N^{\frac{2n}{2n+1}}
\displaystyle\sum_{2k+n|N}
(2k+n)^{-1}
||f||^2_{L^{\tp}(\heiserid)}\,,\\
\end{align*}
whence
\begin{equation}\label{stimatp-2}
||Q_N ||_{\left( L^{\tp}, L^2\right) }\le 
CN^{\frac{n}{2n+1}} \left( d(N)\right)^{\unme}\,.\end{equation}
Thus by applying the Riesz-Thorin interpolation theorem to
$(\ref{stima1-2})$ 
and to
$(\ref{stimatp-2})$ 
we get $(\ref{stimaQNp2})$.
\end{proof}

\begin{remark}
Observe that estimate
($\ref{stimaQNp2}$)
is better than 
($\ref{thang}$) only when $d(N)<1$.

Thus, on the  one hand we are led to seek
arithmetic progressions  $\{N_m\}$ on which
the divisor function $d(N_{m})$, whose behaviour is in general
highly irregular,  is strictly smaller than one. 
On the other one, we are led to inquire about the average size of the norm of 
$Q_N$.
\par\noindent
We remark that, if $n=1$ then $d(N)$ is necessarily greater than one.
\end{remark}

\begin{remark}
Proposition
$\ref{stimaQN}$
reveals the existence of a critical point 
$\tp\in (1,2)$, where the form of the  exponent
 of the eigenvalue $N$
in ($\ref{thangricci}$)
changes.
\end{remark}

In the following  we list some cases in which 
estimate
($\ref{stimaQNp2}$) really
improves  the result in \cite{Thangavelu1}.
First of all,
when  $n\ge 2$ and $N$ is a prime number,
Proposition $\ref{stimaQN}$
 yields the following sharp result.

\begin{proposition}
\label{stimaQNprimi}
Let $n> 2$, $n$ odd.
Let $N$ be a prime number.

\noindent
Then
for every $1\le p\le 2$ 
 \begin{equation}\label{stimaQNpprimi}
||Q_{N}|| _{\left(L^{p}(h^{n}), L^2 (h^{n})
\right) }
\le
\begin{cases}\!
 C \,{N}^{n(\unpi-\unme)-\unme}\,,
&\,\text{if
  $1\le p< \tp$}\cr
C\,{N}^{-\unme(\unpi-\unme)
}
&\,\text{if
  $\tp\le p\le 2$,}\cr
  \end{cases}
 \end{equation}
 with $\tp=
2
{ \frac{2n+1}{2n+3}}$.
Moreover, the above estimate  is sharp.
\end{proposition}
\begin{proof}
($\ref{stimaQNpprimi}$)
follows directly from
($\ref{stimaQNp2}$).

Furthermore,  since in this case
$$||Q_N|| _{\left(L^{p}(h^{n}), L^2 (h^{n})
\right) }\sim
||P_{{1},{\frac{N-n}{2}}}|| _{\left(L^{p}(h^{n}), L^2 (h^{n})
\right) }\,,\;\;\,1\le p\le 2\,,$$
sharpness follows from
Theorem  $\ref{teoremariassuntivo2}$.
\end{proof}

Proposition  $\ref{stimaQNprimi}$ may be generalized to the
case $N=r^{k_0}$, where ${k_0}\in\NN$
and $r$ varies in the set of all prime numbers.

\begin{proposition}
\label{stimaQNpotenzaprimi}
Let $n\ge 2$ be  odd.
Fix a positive integer number $k_0$.
Set $N_r=r^{k_0}$, where 
$r$ varies in the set of all  prime numbers.

\noindent
Then
for every $1\le p\le 2$ 
 \begin{equation}\label{stimaQNpotenzeprimi}
||Q_{N_r}|| _{\left(L^{p}(h^{n}), L^2 (h^{n})
\right) }
\le
\begin{cases}\!
 C \,{N_r}^{n(\unpi-\unme)-\unmeko}\,,
&\,\text{if
  $1\le p< \tp$}\cr
C\,{N_r}^{(n-\unmeko(2n+1))(\unpi-\unme)
}
&\,\text{if
  $\tp\le p\le 2$,}\cr
  \end{cases}
 \end{equation}
 with $\tp=
2
{ \frac{2n+1}{2n+3}}$.
Moreover, 
$(\ref{stimaQNpotenzeprimi})$
is sharp.
\end{proposition}
\begin{proof}
($\ref{stimaQNpotenzeprimi}$)
follows directly from
($\ref{stimaQNp2}$), since
$$d(N_r)=
\frac{1}{r}
+
\frac{1}{r^2}
+\ldots
+
\frac{1}{r^{k_0}}\,\le \frac{2}{r}.$$

To prove that
($\ref{stimaQNpotenzeprimi}$)
is sharp, take the joint eigenfunction $f_0$  for
$L$ and $i^{-1}\partial_t$, with eigenvalues, respectively,
$(2k+n)m=N_r$ and $m=r^{k_0 -1}$, yielding the sharpness
for the joint spectral projection
$P_{r^{k_0 -1}, {\frac{r-n}{2}}}$,
 that is such that
$$|| P_{r^{k_0 -1}, {\frac{q-n}{2}}}||_{(p,2)}\sim
\frac{||f_0||_{p'}}{||f_0||_2}\,.$$
Now we have
\begin{align*}
||Q_N|| _{\left(L^{2}(h^{n}), L^{p'} (h^{n})
\right) }&\ge
\frac{||Q_N f_0|| _{ L^{p'}  }}{|| f_0|| _{ L^{2}  }}=
\frac{|| f_0|| _{ L^{p'}  }}{|| f_0|| _{ L^{2}  }}\sim
|| P_{r^{k_0 -1}, {\frac{r-n}{2}}}||_{(p,2)}\\
&\sim C
r^{n(\unpi-\unme)-\unme}
\left({r^{k_0 -1}}\right)^{n (\unpi-\unme)}\sim 
Cr^{-\unme}{r}^{k_0 n (\unpi-\unme)}\\
&\sim
C{N_r}^{ n (\unpi-\unme)-\unmeko}\,\\
\end{align*}
for all $1\le p\le \tp$.
For $\tp\le p\le 2$ an analogous estimate hold, so that ($\ref{stimaQNpotenzeprimi}$)
is sharp.
\end{proof}

We shall now consider integers of the form
$N_{\ell}:={q_0}^{\ell}$,
where
$q_0$ is a fixed prime number
and $\ell\in\NN$.
The argument of the previous proposition also proves the following.

\begin{proposition}
\label{stimaQNprimofissato}
Let $n= 2$ or $n>2$ odd.
For $n=2$ let $q_0=2$, for $n>2$ let
$q_0$ be a  prime number  strictly greater than  $2$.
Set $N_{\ell}:=
{q_0}^{\ell}$,
$\ell\in\NN$.

\noindent
Then
 \begin{equation}\label{stimaQNfissatoprimo}
||Q_{N_{\ell}}|| _{\left(L^{p}(h^{n}), L^2 (h^{n})
\right) }
\le
 C \,{N_{\ell}}^{n(\unpi-\unme)}\,
\;\,\text{if
  $1\le p\le 2$.}
  \end{equation}
 Moreover, 
$(\ref{stimaQNfissatoprimo})$
is sharp.
\end{proposition}



The above examples show the highly irregular behaviour of $d(N)$, and therefore of
$||Q_N||_{p,2}$.
In order to smooth out
fluctuations 
we introduce appropriate averages of 
joint spectral projectors.
More precisely, we define for $N\in\NN$
\begin{equation}\label{defsommeproiettori}
\Pi_N f:=
\displaystyle\sum_{L=n}^N\,\,
\displaystyle\sum_{(2k+n)|m|=L}P_{m,k} f\,
\end{equation}
and ask what is the behaviour 
of 
$||M_N||_{(p,2)}$,
where
\begin{equation}\label{defmediepr}
M_N f:=\frac{1}{N} \Pi_N
f\,.
\end{equation}

For $p=1$ 
Theorem 
$\ref{teoremariassuntivo2}$ and orthogonality 
yield
\begin{align*}
||\Pi_N f||^2_{L^2(\heisenunouno)}&=
\displaystyle\vert\vert
\displaystyle\sum_{L=n}^N\,\,
\sum_{(2k+n)|m|=L}
P_{m,k}f\vert\vert^2_{L^2(\heisenunouno)}\\
&=
\displaystyle
\sum_{(k,m):\,{(2k+n)|m|}\le N}
\vert\vert
P_{m,k}f\vert\vert^2_{L^2(\heisenunouno)}\\
&\le C\,
\displaystyle
\sum_{(k,m):\,{(2k+n)|m|}\le N}
\left( {2k+n}{}\right)^{n-1
}
|m|^{n}
||f||^2_{L^1(\heisenunouno)}
\\
&\le C\,
\displaystyle
\sum_{\,{m=1}}^N
m^n
\,\sum_{2k+n=n}
^{\left[\frac{N}{m}\right]}
(2k+n)^{n-1}
||f||^2_{L^1(\heisenunouno)}\le
 C\, N^n\cdot N
||f||^2_{L^1(\heisenunouno)}\,,
\end{align*}
whence
\begin{equation}
|| \Pi_N||_{(1,2)}
\le
N^{\frac{n+1}{2}}\,.
\end{equation}
The trivial  $L^2-L^2$ estimate and  Riesz-Thorin interpolation  yield
\begin{equation}
||\Pi_N||_{(p,2)}\label{rettaPi}
\le C\,
N^{{(n+1)}({\unpi}-\unme)}\,
\qquad1\le p\le 2
\end{equation}

Observe that by  using Theorem 
$\ref{teoremariassuntivo2}$
we may obtain
the following estimate
in the critical point $\tp$

\begin{align*}
||\Pi_N f||^2_{L^2(\heisenunouno)}&=
\displaystyle
\sum_{(k,m):\,{(2k+n)|m|}\le N}
\vert\vert
P_{m,k}f\vert\vert^2_{L^2(\heisenunouno)}\\
&\le C\,
\displaystyle
\sum_{(k,m):\,{(2k+n)|m|}\le N}
\left( {2k+n}{}\right)^{2\alpha
}
m^{2\beta}
||f||^2_{L^{\tp}(\heisenunouno)}
\\
&= C\,
\displaystyle
\sum_{\,{m=1}}^N
m^{2\beta}
\,\sum_{2k+n=n}
^{\frac{N}{m}}
(2k+n)^{2\alpha}
||f||^2_{L^{\tp}(\heisenunouno)}
=
 {N}^{2\alpha +1}
\displaystyle
\sum_{\,{m=1}}^N
m^{2\beta-2\alpha-1}
||f||^2_{L^{\tp}(\heisenunouno)}
\\&\le
 C\, N^{2\alpha+2}
||f||^2_{L^{\tp}(\heisenunouno)}\,,
\end{align*}
where we used the fact that
${2\beta-2\alpha}=1$ for all $1\le p\le \tp$,
with $\alpha=\alpha (\unpi, n)$ and $\beta=\beta(\unpi,n)$
given by
($\ref{alfadef}$)
and
($\ref{betadef}$).

Thus 
\begin{equation}\label{tp,2Pi}
||\Pi_N||_{(\tp,2)}
\le C\, N^{\alpha+1}=C\,
N^{\frac{2n+\unme}{2n+1}}\,.
\end{equation}
A comparison between 
($\ref{rettaPi}$)
and
($\ref{tp,2Pi}$)
shows that
in the critical point the estimate given by Riesz-Thorin interpolation is better
than the bound obtained by summing up the estimates for joint spectral projections.

Thus we obtain the following result.

 \begin{proposition}
\label{stimamediahn}
Let $n\ge 1$.
The following $L^p- L^2$ bounds hold for 
$\Pi_N$
and for the average projection operators $M_N$
$$||\Pi_{N}|| _{
\left(
L^{p}(h^{n}), L^2 (h^{n})
\right) }
\le
 C \,{N}^{(n+1)(\unpi-\unme)}\,
\;\,\text{if
  $1\le p\le 2$.}
$$
and  
$$||{}M_{N}|| _{
\left(
L^{p}(h^{n}), L^2 (h^{n})
\right) }
\le
 C \,{N}^{(n+1)(\unpi-\unme)-1}\,
\;\,\text{if
  $1\le p\le 2$.}
  $$
  \end{proposition}

A similar proof also yields the following result about the operators $E_{N_{1}, N_{2}}$,
where
$$E_{N_{1}, N_{2}}:=
\Pi_{ N_{2}}-
\Pi_{N_{1}}\,,\qquad N_1\,,N_2\in \NN\,,\,N_2>N_1\,.$$

 \begin{proposition}
Let $n\ge 1$.
Then
$$||E_{N_{1}, N_{2}}|| _{
\left(
L^{p}(h^{n}), L^2 (h^{n})
\right) }
\le
 C \,
 \left(
  N_{2}^n
  (N_{2}- N_{1})
  \right)^{(\unpi-\unme)}\,
\;\,\text{for all
  $1\le p\le 2$.}
$$
  \end{proposition}
  \begin{remark}
This should be compared to Proposition 3.8
in   \cite{Mueller}.
  \end{remark}

\end{document}